\documentclass[10pt]{article}
\usepackage{amsmath}
\usepackage{graphicx}
\usepackage{amsfonts}
\usepackage{amssymb}
\usepackage{color}

\newcommand{\1}{{\bf 1}}
\newcommand{\la}{\lambda}

\newcommand{\cc}{{\cal C}}

\newcommand{\cf}{{\cal F}}

\newcommand{\al}{\alpha}

\newcommand{\si}{\sigma}

\newcommand{\R}{{\mathbb R}}

\newcommand{\G}{{\mathbb P}}

\newtheorem{theorem}{Theorem}[section]

\newtheorem{lemma}[theorem]{Lemma}

\newtheorem{proposition}[theorem]{Proposition}

\newtheorem{remark}[theorem]{Remark}

\begin{document}
\thispagestyle{empty}

\begin{center}
\huge {

Stochastic Volterra equations
driven by fractional Brownian motion
with Hurst parameter H $>$ 1/2
}

\vspace{.5cm}

\normalsize {\bf Mireia Besal\'u} and {\bf Carles
Rovira$^*$}

{\footnotesize \it

Facultat de Matem\`atiques, Universitat de Barcelona, Gran
Via 585, 08007-Barcelona, Spain.

 {\it E-mail addresses}: mbesalu@ub.edu,
Carles.Rovira@ub.edu}

{$^*$ Corresponding author.}

\end{center}

\begin{abstract}%
In this note we prove an existence and uniqueness result of
solution for stochastic Volterra integral equations driven by a fractional Brownian motion with Hurst
parameter $H > 1/2$, showing also that the solution has finite moments.  The
stochastic integral with respect to the fractional Brownian motion
is a pathwise Riemann-Stieltjes integral.

\end{abstract}


Riemann-Stieltjes integral

{\bf AMS 2000 MSC:} 60H05, 60H07

{\bf Running head:} Stochastic Volterra equations driven by fBm

\setcounter{section}{-1}

\renewcommand{\theequation}{0.\arabic{equation}}
\setcounter{equation}{0}

\section {Introduction}

Consider the stochastic Volterra equation on $\R^d$
\begin{equation}
X(t)= X_0+\int_0^t b(t,s,X(s))ds+\int_0^t\sigma(t,s,X(s))dW_s^H,\quad t\in (0,T],  \label{equation}
\end{equation}
where $W^H=\left\{W^{H,j},j=1,\ldots, m \right\}$ are independent fractional Brownian motions with Hurst parameter $H>\frac{1}{2}$ defined in a complete probability space $(\Omega, \cf, \G)$.

In this paper, the integral with respect W is a pathwise Riemann-Stieltjes integral.
We will define it using a pathwise approach.
Indeed, Young \cite{Y} proved that if we have a stochastic processes $\{u(t), t \ge 0\}$
whose trajectories are $\lambda$-H\"older continuous with $\lambda
> 1-H$, then the Riemann-Stieltjes integral $\int_0^T u(s) dW_s^H$
exists for each trajectory.
Using
the techniques of fractional calculus, Z\"ahle \cite{Z} introduced a generalized Stieltjes integral that coincides with
the Riemann-Stieltjes integral $\int_0^T f dg$ when the functions
$f$ and $g$ are H\"older continuous of orders $\lambda$ and $\beta$,
respectively, with $\lambda + \beta > 1$. Moveover, this generalized Stieltjes integral can be expressed
in terms of a fractional derivative operator.

Following the ideas given by  Z\"ahle \cite{Z} and using the
Riemann-Stieltjes integral, Nualart and Rascanu \cite{N-R} proved a general result on
the existence and uniqueness of solution for a class of multidimensional time dependent stochastic
differential equations driven by a fractional Brownian motion with Hurst parameter $H > 1/2$.
Their proofs begin with a deterministic existence uniqueness theorem based on some a priori estimates
on Lebesgue and  the generalized Stieltjes integral.

In our paper, we extend the results in Nualart and Rascanu
\cite{N-R} to multidimensional stochastic Volterra equations. Using
the Riemann-Stieltjes integral, we give the existence and uniqueness
of a solution to our equation (\ref{equation}). We also show that
the solution has finite moments. Since we follow the methodology
presented in Nualart and Rascanu \cite{N-R}, our main aim is to
obtain precise estimates for Lebesgue and Riemann-Stieltjes Volterra
integrals. Once we have obtained these estimates, the proofs of our
existence and uniqueness results follow exactly the ones in Nualart
and Rascanu \cite{N-R}. Let us note that our results include as a
particular case the results of Nualart and Rascanu.

There are a lot of  references on stochastic differential equations driven by a fractional Bownian motion and many papers about
stochastic Volterra equations (see for instance \cite{AN}, \cite{BM1}, \cite{BM2}, \cite{P}). Nevertheless
the literature about Volterra equations driven by a fractional Brownian
motion is scarce. As far as the authors know, the main references are the papers of Deya and Tindel \cite{DT1}, \cite{DT2}.
In these papers they consider the case $H >1/3$ using an algebraic integration setting. For the case $H>1/2$ they also use a Young integral but they deal with $b=0$ and the set of hypothesis on the coefficients are different and stronger in some sense.

On the other hand, several references have followed the ideas of Nualart and Rascanu for different types of stochastic models (see \cite{FR}, \cite{GN}, \cite{HZ}, \cite{NS}).

The structure of the paper is as follows: in the next section we
state the main result of our paper. In Section 2 we study some
estimates for Lebesgue integrals. Section 3 is devoted to obtain
similar estimations for Riemann-Stieltjes integrals. In Section 4 we
recall the results about deterministic equations and how we apply
them to stochastic equations driven by fractional Brownian motion.
Finally, in Section 5 we give some technical lemmas that we use
throughout the paper.

\setcounter{equation}{0}
\section{Main results}
Let $\frac{1}{2}< H< 1,\,\alpha\in\left(1-H,\frac{1}{2}\right)$. Denote by $W_0^{\alpha,\infty}(0,T;\R^d)$ the space of measurable functions $f:[0,T]\rightarrow\R^d$ such that
\[\left\|f\right\|_{\alpha,\infty}:=\sup_{t\in[0,T]}\left(\left|f(t)\right|+\int_{0}^t\frac{\left|f(t)-f(s)\right|}{(t-s)^{\alpha+1}}ds\right)<\infty.\]
For any $0<\lambda\leq 1$, denote by $C^{\lambda}(0,T;\R^d)$ the space of $\lambda-$H\"older continuous functions $f:[0,T]\rightarrow\R^d$ such that
\[\left\|f\right\|_{\lambda}:=\left\|f\right\|_{\infty}+\sup_{0\leq s<t\leq T}\frac{\left|f(t)-f(s)\right|}{(t-s)^{\lambda}}<\infty,\]
where
\[\left\|f\right\|_{\infty}:=\sup_{s\in[0,T]}\left|f(s)\right|.\]

\vskip 7pt
\noindent
Let us consider the following hypothesis:
\begin{itemize}
\item {\bfseries(H1)} $\sigma:[0,T]^2\times\R^d\rightarrow \R^d\times\R^m$ is a mesurable function such that there exists the derivatives  $\partial_x \sigma(t,s,x)$,  $\partial_t \sigma(t,s,x)$ and $\partial^2_{x,t} \sigma(t,s,x)$ and
there exists some constants  $0<\beta,\,\mu,\, \delta\leq
1$ and for every $N\geq 0$ there exists $K_N>0$ such that the following properties hold:
\begin{enumerate}
\item $\left|\sigma(t,s,x)-\sigma(t,s,y)\right| + \left|\partial_{t}\sigma(t,s,x)-\partial_{t}\sigma(t,s,y)\right|\leq K\left|x-y\right|,$\\$\forall x,y\in\R^d,\;\forall s,t\in[0,T]$,
\item $\left|\partial_{x_i}\sigma(t,s,x)-\partial_{y_i}\sigma(t,s,y)\right|+\left|\partial_{x_i,t}^2\sigma(t,s,x)-\partial_{y_i,t}^2\sigma(t,s,y)\right|\leq K_N\left|x\right.\\ \left.-y\right|^\delta,\forall |x|,|y|\leq N,\;\forall s,t\in[0,T],\; i=1\ldots d$,
\item $\left|\sigma(t_1,s,x)-\sigma(t_2,s,x)\right|+\left|\partial_{x_i}\sigma(t_1,s,x)-\partial_{x_i}\sigma(t_2,s,x)\right|
\leq K\left|t_1-\right. \\ \left.t_2\right|^\mu,$ $\forall
x\in\R^d,\;\forall t_1,t_2,s\in[0,T],\; i=1\ldots d$,
\item $\left|\sigma(t,s_1,x)-\sigma(t,s_2,x)\right|+\left|\partial_{t}\sigma(t,s_1,x)-\partial_{t}\sigma(t,s_2,x)\right|\leq K\left|s_1-s_2\right|^\beta,\\ \forall x\in\R^d,\;\forall s_1,s_2,t\in[0,T]$,
\item $\left|\partial_{x_i,t}^2\sigma(t,s_1,x)-\partial_{x_i,t}^2\sigma(t,s_2,x)\right|+\left|\partial_{x_i}\sigma(t,s_1,x)-\partial_{x_i}\sigma(t,s_2,x)\right|\\ \leq K\left|s_1-s_2\right|^\beta$,$\qquad\forall x\in\R^d,\;\forall s_1,s_2,t\in[0,T],\, i=1\ldots d$.
\end{enumerate}
\item {\bfseries(H2)} $b:[0,T]^2\times \R^d \rightarrow \R^d$ is a measurable function such that there exists $b_0\in L^\rho([0,T]^2;\R^d)$ with $\rho\geq 2$, $0<\mu\leq 1$ and $\forall N\geq 0$ there exists $L_N>0$ such that:
\begin{enumerate}
\item $\left|b(t,s,x)-b(t,s,y)\right|\leq L_N\left|x-y\right|,\; \forall |x|,|y|\leq N,\,\forall s,t\in[0,T]$,
\item $\left|b(t_1,s,x)-b(t_2,s,x)\right|\leq L\left|t_1-t_2\right|^\mu,\; \forall x\in\R^d,\,\forall s,t_1,t_2\in[0,T]$,
\item $\left|b(t,s,x)\right|\leq L_0|x|+b_0(t,s),\quad \forall x\in\R^d,\;\forall s,t\in[0,T]$,
\item $\left|b(t_1,s,x_1)-b(t_1,s,x_2)-b(t_2,s,x_1)+b(t_2,s,x_2)\right|\leq L_N |t_1-t_2||x_1-x_2|,\quad
\forall |x_1|,|x_2|\leq N,\;\forall t_1,t_2,s\in[0,T]$.
\end{enumerate}
\item {\bfseries(H3)} There exists $\gamma\in[0,1]$ and $K_0>0$ such that
\[|\si(t,s,x)|\leq K_0(1+|x|^\gamma), \;\forall x\in\R^d,\,\forall s,t\in[0,T].\]
\end{itemize}
\vskip 6pt
\begin{remark}
Actually, we can consider $\sigma$ and $b$ defined only in the set
$D \times \R^d$ with $D=\{(t,s)\in [0,T]^2; s \le t \}$.
\end{remark}
\vskip 6pt
Under these assumptions we are able to prove that our problem admits a unique solution. The result of existence and uniqueness reads as follows:
\vskip 6pt
\begin{theorem}
Assume that $X_0$ is a $\R^d$-valued random variable and that $b$ and $\sigma$ satisfy hypothesis
{\upshape\bfseries (H1)} and {\upshape\bfseries (H2)} respectively with
 $\beta>1-H$, $\delta>\frac{1}{H}-1$, $ \min \{\beta, \frac{\delta}{1+\delta} \} >1- \mu.$ Set $\alpha_0:=\min\left\{\frac{1}{2}, \beta,
  \frac{\delta}{1+\delta}\right\}$.\\
Then if $\alpha\in((1-H) \vee (1-\mu),\alpha_0)$ and
$\rho\leq \frac{1}{\alpha}$, equation (\ref{equation}) has an unique
solution
\[X\in L^0(\Omega,\cf,\G;W_0^{\alpha,\infty}(0,T;\R^d))\]
and for $P-almost\;all\,\omega\in\Omega,\;X(\omega,\cdot)\in C^{1-\alpha}(0,T;\R^d)$.

Moreover, if $\alpha \in ((1-H) \vee (1-\mu), \al_0 \vee
(2-\gamma)/4),\; \rho \ge 1/\alpha,\; X_0 \in
L^\infty(\Omega,\cf,\G;\R^d)$ and {\bf (H3)} holds then $E ( \Vert X
\Vert^p_{\al,\infty} ) < \infty,\; \forall p \ge 1.$
\end{theorem}

\vskip 6pt

 In order to prove these results, we  need to introduce a
new norm in the space $W_0^{\alpha,\infty}(0,T;\R^d)$ that is for
any $\lambda\geq 1$:
\[\left\|f\right\|_{\alpha,\lambda}:=\sup_{t\in[0,T]}\exp(-\lambda t)\left(\left|f(t)\right|+\int_{0}^t \frac{\left|f(t)-f(s)\right|}{(t-s)^{\alpha+1}}ds\right).\]
It is easy to check that, for any $\lambda\geq 1$, this norm is
equivalent to $\left\|f\right\|_{\alpha,\infty}$.

\renewcommand{\theequation}{2.\arabic{equation}}
\setcounter{equation}{0}
\section{Lebesgue integral}
Let us consider first the ordinary Lesbesgue integral. Given $f:[0,T]^2\rightarrow\R^d$ a measurable function we define
\[F_t(f)=\int_0^t f(t,s)ds. \]
\begin{proposition}
Let $0<\al<\frac{1}{2}$ and $f:[0,T]^2\rightarrow\R^d$ a measurable
function such that for all $0<s \le t_1, t_2,$
$|f(t_1,s)-f(t_2,s)|\leq L|t_1-t_2|^\mu \,$ with $\mu
> \alpha.$
If ${\displaystyle\sup_{t\in[0,T]}\int_0^t\frac{|f(t,s)|}{(t-s)^\al}ds<\infty}$ then $F_\cdot(f)\in W_0^{\al,\infty}(0,T;\R^d)$ and
\begin{equation}
\left|F_t(f)\right|+\int_0^t\frac{\left|F_t(f)-F_s(f)\right|}{(t-s)^{\al+1}}ds\leq
C_{\al,T}^{(1)}
\int_0^t\frac{|f(t,s)|}{(t-s)^\al}ds+C_{\al,L,\mu}^{(2)}t^{
1+\mu-\al} \label{cotaf1}.
\end{equation}
\end{proposition}
{\bf Proof:}
We have that
\begin{eqnarray*}
\left|F_t(f)\right|&+&\int_0^t\frac{|F_t(f)-F_s(f)|}{(t-s)^{\al+1}}ds \le \left|\int_0^t f(t,u)du\right|\\
&&+\int_0^t\frac{\int_0^s \left|f(t,u)-f(s,u)\right|du}{(t-s)^{\al+1}}ds+
\int_0^t \frac{\int_s^t |f(t,u)| du}{(t-s)^{\al+1}}ds\\&\leq& T^\al\int_0^t\frac{\left|f(t,u)\right|}{(t-u)^\al}du
+L\int_0^t\frac{s}{(t-s)^{1-\mu+\al}}ds\\&&+\int_0^t\int_0^u\frac{|f(t,u)|}{(t-s)^{\al+1}}dsdu\\
&\leq&
\left(T^\al+\frac{1}{\al}\right)\int_0^t\frac{|f(t,u)|}{(t-u)^\al}du+\frac{L}{(\mu-\al)}
t^{1+ \mu-\al},
\end{eqnarray*}
so (\ref{cotaf1}) holds with
$C_{\al,T}^{(1)}=\left(T^\al+\frac{1}{\al}\right)$ and
$C_{\al,L,\mu}^{(2)} =\frac{L}{(\mu-\al)}$ and we also have proved
that $F_.(f)\in W_0^{\al,\infty}(0,T;\R^d)$. \hfill $\Box$ \vskip
5pt \noindent Given $f:[0,T] \to \R^d$ let us now define
$$F_t^{(b)}(f)=\int_0^t b(t,s,f(s))ds.$$

\begin{proposition} \label{Prop4.4}
Assume that b satisfies {\bfseries\upshape (H2)} with
$\rho=\frac{1}{\al}$ and $\mu >(1-\alpha) \vee
\alpha$.
\begin{enumerate}
\item If $f\in W_0^{\al,\infty}(0,T;\R^d)$ then $F_\cdot^{(b)}(f)\in \cc^{1-\al}(0,T;\R^d)$ and
\begin{eqnarray}
\left\|F^{(b)}(f)\right\|_{1-\al}&\leq& d^{(1)}\left(1+\left\|f\right\|_\infty\right),\label{norma1} \\
\left\|F^{(b)}(f)\right\|_{\al,\la}&\leq& \frac{d^{(2)}}{\la^{1-2\al}}\left(1+\left\|f\right\|_{\al,\la}\right)\label{norma2},
\end{eqnarray}
for all $\la\geq 1$, where $d^{(1)}$ and $d^{(2)}$ are positive constants depending only on $\mu, \al, T, L, L_0$ and a constant $B_{0,\al}$ that depends on $b$.
\item If $f,h\in W_0^{\al,\infty}\left(0,T;\R^d\right)$ are such that $\left\|f\right\|_\infty\leq N,\;\left\|h\right\|_\infty\leq N,$ then
\[\left\|F^{(b)}(f)-F^{(b)}(h)\right\|_{\al,\la}\leq \frac{d_N}{\la^{1-\al}}\left\|f-h\right\|_{\al,\la},\]
for all $\la\geq 1$, where $d_N$ depends on $\al,T$ and $L_N$ from {\bfseries\upshape (H2)}.
\end{enumerate}
\end{proposition}
{\bf Proof:}
In order to simplify the presentation we will assume  $d=1$. Let $f\in W_0^{\al,\infty}(0,T)$, then for $0\leq s<t\leq T$
\begin{eqnarray}
\left|F_t^{(b)}(f)-F_s^{(b)}(f)\right|&\leq& \int_0^s
\left|b(t,u,f(u))-b(s,u,f(u))\right|du\nonumber\\&&+ \int_s^t
|b(t,u,f(u))|du\nonumber\\&\leq& Ls(t-s)^{
\mu}+\int_s^t\left(L_0|f(u)|+b_0(t,u)\right)du\nonumber\\&\leq&
(t-s)^{1-\al}\left(LT^{ \mu+\al}+L_0T^\al\left\|f\right\|_\infty+
B_{0,\al}\right) \label{incb1},
\end{eqnarray}
where $B_{0,\al}:= \sup_{t \in [0,T]} \left( \int_0^t \vert b_0 (t,u) \vert^{1/\al} du \right)^\alpha$.
The same computations  for $s=0$ give us
\[\left|F_t^{(b)}(f)\right|\leq L_0t\left\|f\right\|_\infty+B_{0,\al}t^{1-\al}.\]
Hence
\begin{equation*}
\left\|F^{(b)}(f)\right\|_{1-\al}\leq
L_0(T+T^\al)\left\|f\right\|_\infty+B_{0,\al}(1+T^{1-\al})+LT^{
\mu+\al},
\end{equation*}
and (\ref{norma1}) holds with $d^{(1)}=(1+T^{1-\al})(B_{0,\al}+T^\al L_0)+LT^{ \mu+\al}$.\\
By (\ref{cotaf1}) we have
\begin{eqnarray}
\left|F^{(b)}_t(f)\right|&+&\int_0^t\frac{\left|F_t^{(b)}(f)-F_s^{(b)}(f)\right|}{(t-s)^{\al+1}}ds\nonumber\\
&\leq&
C_{\al,T}^{(1)}\int_0^t\frac{|b(t,s,f(s))|}{(t-s)^\al}ds+C_{\al,L}^{(2)}t^{
1+\mu-\al}\nonumber\\&\leq&
C_{\al,T}^{(1)}\int_0^t\frac{L_0|f(s)|+b_0(t,s)}{(t-s)^\al}ds+C_{\al,L}^{(2)}t^{
1+\mu-\al}\nonumber\\&\leq&
C_{\al,T}^{(1)}\left(L_0\int_0^t\frac{|f(s)|}{(t-s)^\al}ds+
B_{0,\al}\left(\frac{1-\al}{1-2\al}\right)^{1-\al}t^{1-2\al}\right)\nonumber\\&&
+C_{\al,L}^{(2)}t^{ 1+\mu-\al}. \label{cotamom}
\end{eqnarray}
Then using that
\begin{equation}
\int_0^t\frac{e^{-\la(t-s)}}{(t-s)^\al}ds\leq
\la^{\al-1}\Gamma(1-\al)\;\quad\text{and}\; \quad
\sup_{t\in[0,T]}t^\mu e^{-\la t}\leq
\left(\frac{\mu}{\la}\right)^\mu e^{-\mu},
\label{sempre}\end{equation}
 we get
that
\begin{eqnarray*}
\left\|F^{(b)}(f)\right\|_{\al,\la}&\leq&
C_{\al,T}^{(1)}L_0\left\|f\right\|_{\al,\la}\sup_{t\in[0,T]}
\int_0^t\frac{e^{-\la(t-s)}}{(t-s)^\al}ds\\&&+C_{\al,L}^{(2)}\sup_{t\in[0,T]}e^{-\la
t}t^{ 1+\mu-\al}
\\&&+C_{\al,T}^{(1)}B_{0,\al}\left(\frac{1-\al}{1-2\al}\right)^{1-\al}\sup_{t\in[0,T]}e^{-\la t}t^{1-2\al}\\
&\leq&
C_{\al,T}^{(1)}L_0\Gamma(1-\al)\la^{\al-1}\left\|f\right\|_{\al,\la}\\&&+C_{\al,L}^{(2)}e^{\al-\mu-1}(
1+\mu-\al)^{ 1+\mu-\al} \la^{\al-
1-\mu}\\&&+C_{\al,T}^{(1)}B_{0,\al}\frac{(1-\al)^{1-\al}}{(1-2\al)^\al}e^{2\al-1}\la^{2\al-1}\\&\leq&
d^{(2)}\la^{2\al-1}\left(1+\left\|f\right\|_{\al,\la}\right).
\end{eqnarray*}
So (\ref{norma2}) holds with
\[d^{(2)}=C_{\al,T}^{(1)}L_0\Gamma(1-\al)+C_{\al,T}^{(1)}B_{0,\al}\frac{(1-\al)^{1-\al}}{(1-2\al)^\al}e^{2\al-1}
+C_{\al,L}^{(2)}e^{\alpha-\mu-1}(1+\mu-\al)^{
1+\mu-\al}.\] Now, consider $f,h\in W_0^{\al,\infty}(0,T)$ such
that $\left\|f\right\|_\infty\leq N,\;\left\|h\right\|_\infty\leq
N$. We obtain that
\begin{eqnarray}
\left|F^{(b)}_t(f)-F^{(b)}_t(h)\right|&\leq& \int_0^t \left|b(t,u,f(u))-b(t,u,h(u))\right|du \nonumber\\&\leq& L_N\int_0^t \left|f(u)-h(u)\right|du \label{incb2}
\end{eqnarray}
and
\begin{equation}\begin{array}{l}
\displaystyle
|F^{(b)}_t(f)-F^{(b)}_t(h)-F^{(b)}_s(f)+F^{(b)}_s(h)|\leq\int_s^t\left|b(t,u,f(u))-b(t,u,h(u))\right|du\\ [4mm]\displaystyle\qquad\quad+\int_0^s \left|b(t,u,f(u))-b(t,u,h(u))-b(s,u,f(u))+b(s,u,h(u))\right|du\\ [3mm]\displaystyle\qquad\quad\leq L_N\int_s^t |f(u)-h(u)|du+L_N |t-s|\int_0^s |f(u)-h(u)|du.
\end{array}\label{incb3}\end{equation}

Then, using (\ref{incb2}) and (\ref{incb3}) we have

\[\begin{array}{l}
\displaystyle \left|F^{(b)}_t(f)-F^{(b)}_t(h)\right|+\int_0^t\frac{|F^{(b)}_t(f)-F^{(b)}_t(h)-F^{(b)}_s(f)+F^{(b)}_s(h)|}{(t-s)^{\al+1}}ds\\[4mm]\displaystyle\qquad\qquad \leq L_N\int_0^t \left|f(u)-h(u)\right|du+L_N \int_0^t \frac{\int_0^s |f(u)-h(u)| du}{(t-s)^\al}ds\\[4mm]\displaystyle\qquad\qquad\quad +L_N\int_0^t\frac{\int_s^t |f(u)-h(u)| du }{(t-s)^{\al+1}}ds\\ [4mm]\displaystyle\qquad\qquad \leq L_N\int_0^t \left|f(u)-h(u)\right|du+\frac{L_N}{1-\al}\int_0^t\frac{|f(u)-h(u)|}{(t-u)^{\al-1}}du\\ [4mm]\displaystyle\qquad\qquad\quad+\frac{L_N}{\al}\int_0^t\frac{|f(u)-h(u)|}{(t-u)^{\al}}du.
\end{array}\]
Finally,
\begin{eqnarray*}
\left\|F^{(b)}(f)-F^{(b)}(h)\right\|_{\al,\la}&\leq&\sup_{t\in[0,T]}\left(L_N \left( 1+ \frac{T^{1-\alpha}}{1-\alpha} \right) \int_0^t e^{-\la(t-u)}du
\right.\\&&+\left.\frac{L_N}{\al}\int_0^t\frac{e^{-\la(t-u)}}{(t-u)^{\al}}du\right)\left\|f-h\right\|_{\al,\la}\\&\leq&
L_N       \left (\frac{1}{\la} \left( 1+ \frac{T^{1-\alpha}}{1-\alpha} \right) +\frac{\Gamma(1-\al)}{\al\la^{1-\al}}\right)\left\|f-h\right\|_{\al,\la}\\&\leq& \frac{d_N}{\la^{1-\al}}\left\|f-h\right\|_{\al,\la},
\end{eqnarray*}
where ${\displaystyle d_N= L_N \left( 1+ \frac{T^{1-\alpha}}{1-\alpha} + \frac{\Gamma(1-\al)}{\al} \right).}$
\hfill $\Box$
\renewcommand{\theequation}{3.\arabic{equation}}
\setcounter{equation}{0}
\section{Riemman-Stieltjes integral}
Let us now consider the Riemman-Stieltjes integral introduced by
Z\"ahle, which is based on fractional integrals and derivatives. We
refer the reader to Z\"ahle \cite{Z} and the references therein for
a detailed account about this generalized integral and the
relationships with fractional calculus. Here we just recall some
basic results. \vskip 7pt Fix a parameter $0<\al<\frac{1}{2}$.
Denote by $W_T^{1-\al,\infty}(0,T)$ the space of measurable
functions $g:[0,T]\rightarrow\R$ such that
\[\left\|g\right\|_{1-\al,\infty,T}:=\sup_{0<s<t<T}\left(\frac{\left|g(t)-g(s)\right|}{(t-s)^{1-\al}}+\int_s^t\frac{\left|g(y)-g(s)\right|}{(y-s)^{2-\al}}dy\right)<\infty.\]
Moreover if g belongs to $W_T^{1-\al,\infty}(0,T)$ we define
\begin{eqnarray*}
\Lambda_\al(g)&:=&\frac{1}{\Gamma(1-\al)}\sup_{0<s<t<T}\left|\left(D^{1-\al}_{t^-}
g_{t^-}\right)(s) \right|\\&\leq&
\frac{1}{\Gamma(1-\al)\Gamma(\al)}\left\|g\right\|_{1-\al,\infty,T}<\infty,
\end{eqnarray*}
where $(D^{1-\al}_{t^-} g_{t^-})(s)$ is a Weyl derivative.  We also
denote by $W_0^{\al,1}(0,T)$ the space of measurable functions f on
$[0,T]$ such that,
\[\left\|f\right\|_{\al,1}:=\int_0^T\frac{\left|f(s)\right|}{s^\al}ds+\int_0^T\int_0^s \frac{\left|f(s)-f(y)\right|}{(s-y)^{\al+1}}dyds<\infty.\]
Then, if $f$ is a function in the space $W_0^{\al,1}(0,T)$ and $g$
belongs to $W_T^{1-\al,\infty}(0,T)$,  the integral $\int_0^t
f(s)dg_s$ exists for all $t\in[0,T]$ and we can define
\[\int_0^t f(s)dg_s=\int_0^T f(s) \1_{(0,t)}(s)dg_s.\]
Furthermore the following estimate holds
\[\left|\int_0^t f(s)dg_s\right|\leq \Lambda_\al (g)\left\|f\right\|_{\al,1}. \]

Given  $f:[0,T]^2\rightarrow\R$  such that for any $t\in [0,T]$, $f(t,\cdot) \in  W_0^{\al,1}([0,T])$
we can also consider the integral
\[G_t(f)=\int_0^t f(t,s)dg_s=\int_0^T f(t,s) \1_{(0,t)}(s)dg_s ,\]
and the estimate
\[\left|\int_0^t f(t,s)dg_s\right|\leq \Lambda_\al (g)\left\|f(t,\cdot)\right\|_{\al,1}. \]
\begin{proposition}
Let $g\in W_T^{1-\al,\infty}(0,T)$ and $f:[0,T]^2\rightarrow\R^d$  such that
 $f(t,\cdot)\in W_0^{\al,1}(0,T)$ for all $t\in[0,T]$ and such that $|f(t_1,s)-f(t_2,s)|\leq K(s)|t_1-t_2|^{
\mu},$ with ${\mu>\alpha}$. Then for all $s<t$, the following
estimates hold
\begin{eqnarray}
\left|G_t(f)-G_s(f)\right|&\leq&
\Lambda_\al(g)\left(|t-s|^{
\mu}\int_0^s\frac{K(u)}{u^\al}du+\int_s^t\frac{\left|f(t,u)\right|}{(u-s)^\al}du\right.\nonumber\\&&\left.+\al\int_0^s\int_0^u\frac{|f(t,u)-f(s,u)-f(t,y)+f(s,y)|}{(u-y)^{\al+1}}dydu\right.\nonumber\\&&\left.+\al\int_s^t\int_s^u\frac{\left|f(t,u)-f(t,y)\right|}{(u-y)^{\al+1}}dydu\right)
\label{cotasi1}
\end{eqnarray}
and \begin{equation}\begin{array}{l} \displaystyle
\left|G_t(f)\right|+\int_0^t\frac{\left|G_t(f)-G_s(f)\right|}{(t-s)^{\al+1}}ds
\leq \Lambda_\al(g)\left(C_{\al}^{(3)}
\int_0^t\frac{K(u)}{u^\al}(t-u)^{{
\mu}-\al}du\right.\\[4mm]
\displaystyle\qquad+\left. C_{\al,T}^{(4)}\int_0^t\left((t-u)^{-2\al}+u^{-\al}\right)\left(|f(t,u)|+\int_0^u\frac{|f(t,u)-f(t,y)|}{(u-y)^{\al+1}}dy\right)du\right.\\[4mm]
\displaystyle\qquad+\left.\al\int_0^t\int_0^s\int_0^u\frac{\left|f(t,u)-f(s,u)-f(t,y)+f(s,y)\right|}{(u-y)^{\al+1}(t-s)^{\al+1}}dyduds\right).
\end{array}\label{cotasi2}\end{equation}
%
\end{proposition}
{\bf Proof:}
We have
\begin{eqnarray}
\left|G_t(f)-G_s(f)\right|&\leq& \left|\int_0^s\left(f(t,u)-f(s,u)\right)dg_u\right|+\left|\int_s^t f(t,u)dg_u\right|\nonumber\\
&\leq&\Lambda_\al(g)\left(\int_0^s\frac{\left|f(t,u)-f(s,u)\right|}{u^\al}du+\int_s^t\frac{\left|f(t,u)\right|}{(u-s)^\al}du\right.\nonumber\\
&&\left.+\al\int_0^s\int_0^u\frac{|f(t,u)-f(s,u)-f(t,y)+f(s,y)|}{(u-y)^{\al+1}}dydu\right.\nonumber\\
&&\left.+\al\int_s^t\int_s^u\frac{\left|f(t,u)-f(t,y)\right|}{(u-y)^{\al+1}}dydu\right)\nonumber\\
&\leq& \Lambda_\al(g)\left(|t-s|^{
\mu}\int_0^s\frac{K(u)}{u^\al}du+\int_s^t\frac{\left|f(t,u)\right|}{(u-s)^\al}du\right.\nonumber\\
&&\left.+\al\int_0^s\int_0^u\frac{|f(t,u)-f(s,u)-f(t,y)+f(s,y)|}{(u-y)^{\al+1}}dydu\right.\nonumber\\
&&\left.+\al\int_s^t\int_s^u\frac{\left|f(t,u)-f(t,y)\right|}{(u-y)^{\al+1}}dydu\right).\label{eq0}
\end{eqnarray}
So (\ref{cotasi1}) holds.
\noindent
Moreover
\begin{equation}\begin{array}{l}
\displaystyle
\int_0^t\frac{\left|G_t(f)-G_s(f)\right|}{(t-s)^{\al+1}}ds\leq
\Lambda_\al(g)\left(\int_0^t\int_0^s
\frac{K(u)}{(t-s)^{\al{
+1-\mu}}u^{\al}}duds\right.\\[4mm]
\displaystyle\quad +\al\int_0^t\int_0^s\int_0^u\frac{\left|f(t,u)-f(s,u)-f(t,y)+f(s,y)\right|}{(u-y)^{\al+1}(t-s)^{\al+1}}dyduds \\[4mm]
\displaystyle\quad+\left.\int_0^t (t-s)^{-\al-1}\left(\int_s^t\frac{|f(t,u)|}{(u-s)^\al}du+\al\int_s^t\int_s^u\frac{|f(t,u)-f(t,y)|}{(u-y)^{\al+1}}dydu\right)ds\right),
\end{array}\label{eq2}\end{equation}
\noindent
where the first term can be written as
\begin{eqnarray}
\int_0^t\int_0^s\frac{K(u)}{u^\al(t-s)^{\al{
+1-\mu}}}duds&=&\int_0^t\int_u^t \frac{K(u)}{u^\al(t-s)^{\al
{+ 1-\mu}}}dsdu\nonumber\\&=&\frac{1}{{
\mu}-\al}\int_0^t \frac{K(u)}{u^\al(t-u)^{\al-{ \mu}}}du.
\label{eq3}
\end{eqnarray}
Notice also that using the same computations as in  Proposition 4.1 in \cite{N-R} we obtain that
\begin{equation}\begin{array}{l}
\displaystyle \int_0^t (t-s)^{-\al-1}\left(\int_s^t\frac{|f(t,u)|}{(u-s)^\al}du+\al\int_s^t\int_s^u
\frac{|f(t,u)-f(t,y)|}{(u-y)^{\al+1}}dydu\right)ds\label{eq4}\\[4mm]
\displaystyle\quad \leq
B(2\al,1-\al)\int_0^t\frac{\left|f(t,u)\right|}{(t-u)^{2\al}}du+\int_0^t\int_0^u
\frac{|f(t,u)-f(t,y)|}{(u-y)^{\al+1}}(t-y)^{-\al}dydu,
\end{array}\end{equation}
where $B(2 \alpha, 1 - \alpha)$ is the Beta function. Let us recall
that
 \begin{equation}
 B(p,q)=\int_0^1
t^{p-1}(1-t)^{q-1}dt=\int_0^\infty\frac{t^{p-1}}{(1+t)^{p+q}}dt.\label{semprebeta}
\end{equation}
\vskip 5pt \noindent Finally, the case  (\ref{eq0}) when $s=0$ gives
us
\begin{equation}
\left|G_t(f)\right|\leq \Lambda_\al(g)\left(\int_0^t\frac{\left|f(t,u)\right|}{u^\al}du+\al\int_0^t\int_0^u\frac{\left|f(t,u)-f(t,y)\right|}{(u-y)^{\al+1}}dydu\right). \label{eq1}
\end{equation}
\vskip 5pt
\noindent
So using the estimates (\ref{eq1}),(\ref{eq2}),(\ref{eq3}) and(\ref{eq4}) the inequality (\ref{cotasi2}) holds
 with $C_{\al}^{(3)}=\frac{1}{{
\mu}-\al}$ and $C_{\al,T}^{(4)}=\max(B(2\al,1-\al),1)+T^{\al}$.
\hfill $\Box$

\vskip 8pt
\noindent Given $f:[0,T] \to \R^d$ let us define now
$$G_t^{(\sigma)} (f) = \int_0^t \sigma(t,s,f(s)) dg_s.$$

\vskip 7pt
\begin{proposition} \label{Prop4.2}
Let $g\in W_T^{1-\al,\infty}(0,T)$. Assume that $\sigma$
satisfies {\bfseries\upshape (H1)} with $\beta>\alpha > 1-\mu$.
\begin{enumerate}
\item If $f\in W_0^{\al,\infty}(0,T;\R^d)$ then $$G^{(\sigma)}(f)\in \cc^{1-\al}(0,T;\R^d)\subseteq W_0^{\al,\infty}(0,T;\R^d).$$ Moreover,
\begin{eqnarray}
\left\|G^{(\sigma)}(f)\right\|_{1-\al}&\leq& \Lambda_\al(g)d^{(3)}\left(1+\left\|f\right\|_{\al,\infty}\right),\label{cotasi3}\\
\left\|G^{(\sigma)}(f)\right\|_{\al,\la}&\leq& \frac{\Lambda_\al(g)d^{(4)}}{\la^{1-2\al}}\left(1+\left\|f\right\|_{\al,\la}\right),\label{cotasi4}
\end{eqnarray}
for all $\la\geq 1$ where the constants $d^{(i)}$ for $i=3,4$ depend
only on $\al,\;\beta,\\{ \mu,}\; K,T$ and $N$.
\item If $f,h\in W_0^{\al,\infty}(0,T;\R^d)$ are such that $\left\|f\right\|_\infty\leq N,\, \left\|h\right\|_\infty\leq N$, then
\begin{equation}
\left\|G^{(\sigma)}(f)-G^{(\sigma)}(h)\right\|_{\al,\la}\leq \frac{\Lambda_\al(g)d'_N}{\la^{1-2\al}}(1+\Delta(f)+\Delta(h))\left\|f-h\right\|_{\al,\la}, \label{cotasi5}
\end{equation}
for all $\la\geq 1$, where
\begin{equation*}
\Delta(f)=\sup_{u\in[0,T]}\int_0^u\frac{|f(u)-f(s)|^\delta}{(u-s)^{\al+1}}ds,
\end{equation*}
and the constant $d'_N$ depends only on $\al,\,\beta, {
\mu,}\,N,\,K$ and $T$.
\end{enumerate}
\end{proposition}
{\bf Proof:}
We will assume $d=m=1$ in order to simplify the presentation of the proof. First we see that if $f\in W_0^{\al,\infty}(0,T)$ then $\si (t,\cdot,f(\cdot)) \in W_0^{\al,\infty}([0,T])$ for all $t \in [0,T]$. Indeed:
\begin{eqnarray*}
\left|\si(t,r,f(r))\right|&+&\int_0^r\frac{\left|\si(t,r,f(r))-\si(t,s,f(s))\right|}{(r-s)^{\al+1}}ds\\
&\leq& K(t^{
\mu}+r^\beta+|f(r)|)+|\si(0,0,0)|+K\int_0^r\frac{|f(r)-f(s)|}{(r-s)^{\al+1}}ds\\&&+K\int_0^r
(r-s)^{\beta-\al-1}ds\\&\leq& K\left(t^{
\mu}+r^\beta+\frac{r^{\beta-\al}}{\beta-\al}\right)+|\si(0,0,0)|\\&&+K\left(|f(r)|+\int_0^r\frac{|f(r)-f(s)|}{(r-s)^{\al+1}}ds\right).
\end{eqnarray*}
So, for all $t$
\begin{equation}\left\|\si(t,\cdot,f(\cdot))\right\|_{\al,\infty}\leq K^{(2)}+K\left\|f\right\|_{\al,\infty},\label{novaho}\end{equation}
with ${\displaystyle K^{(2)}=K\left(T^{
\mu}+T^\beta+\frac{T^{\beta-\al}}{\beta-\al}\right)+|\si(0,0,0)|}$.
\vskip 7pt \noindent Now if $f\in W_0^{\al,\infty}(0,T)$ under
assumptions {\bfseries\upshape (H1)}  from (\ref{eq1}) we have,
\begin{eqnarray*}
\left\|G^{(\sigma)} (f)\right\|_\infty&\leq& \sup_{t\in[0,T]}\Lambda_\al(g)\left(\int_0^t\frac{|\sigma(t,u,f(u))|}{u^\al}du\right.\\&&\left.+\al\int_0^t\int_0^u\frac{|\sigma(t,u,f(u))-\sigma(t,y,f(y))|}{(u-y)^{\al+1}}dydu\right)\\
&\leq& \Lambda_\al(g) \left(\frac{T^{1-\al}}{1-\al}+\al T\right)\sup_{t \in [0,T]} \left\|\sigma(t,\cdot,f(\cdot))\right\|_{\al,\infty}.
\end{eqnarray*}
If we come back to (\ref{cotasi1}) with $K(u)=K$ and using Lemma \ref{lema2} we have,
\[\begin{array}{l}
\displaystyle
\left|G_t^{(\sigma)}(f)-G_s^{(\sigma)}(f)\right|\\[4mm]
\displaystyle\qquad\leq\Lambda_\al(g)\left( \vert t- s
\vert^{ \mu} K \int_0^t u^{-\alpha} du + \int_s^t
\frac{\Vert \sigma (t,\cdot,f(\cdot) ) \Vert_\infty }{(u-s)^\alpha}
du \right.
\\[4mm]
\displaystyle\qquad\qquad
\left. + \alpha \int_0^s \int_0^u K \frac{\vert t-s \vert \left( \vert u -y \vert^\beta + \vert f(u)-f(y) \vert \right)}{(u-y)^{\alpha+1}} dy du
\right.
\\[4mm]
\displaystyle\qquad\qquad
\left. + \alpha \int_s^t \int_s^u K \frac{ \vert u -y \vert^\beta + \vert f(u)-f(y)\vert }{(u-y)^{\alpha+1}} dy du
\right.
\\[4mm]
\displaystyle\qquad\leq (t-s)^{1-\al}\Lambda_\al(g) K
\left(\frac{{ T^{\mu}}}{1-\alpha} + \frac{\Vert
\sigma(t,\cdot,f(\cdot) \Vert_\infty}{1-\al} \right.
\\[4mm]
\displaystyle\qquad\qquad \left. + \alpha T^{1+\beta} (
\frac{1}{(\beta-\alpha)(1+\beta-\alpha)} + \Vert f
\Vert_{\alpha,\infty} ) \right.
\\[4mm]
\displaystyle\qquad\qquad \left. + \alpha 
 ( \frac{T^{\beta}}{(\beta-\alpha)} +T^{\alpha}  \Vert f \Vert_{\alpha,\infty} )\right)\\[4mm]
\displaystyle\qquad\leq (t-s)^{1-\al}\Lambda_\al(g)
K_{\al,T}^{(1)}\left(1+\Vert \sigma(t,\cdot,f(\cdot)
\Vert_{\alpha,\infty} + \Vert f \Vert_{\alpha,\infty}\right),
\end{array}\]
where $$K_{\al,T}^{(1)}=5K \left( \frac{1}{1+\al}+\frac{1}{\beta-\al} +\frac{1}{(\beta-\al)(1+\beta-\alpha)} \right)
 (1 +
\alpha T^{1+\beta} {}).$$ So, using (\ref{novaho})
we can  deduce that $G^{(\si)}(f)\in \cc^{1-\al}(0,T)$, and
(\ref{cotasi3}) holds.

\vskip 7pt \noindent Let us study now the norm $\Vert \cdot
\Vert_{\alpha,\lambda}$. Set
$$
\Sigma (t,s,u,y,f):=\si(t,u, f(u))-\si(s,u, f(u))-\si(t,y,
f(y))+\si(s,y,f(y)).
$$
Using (\ref{cotasi2}) with $K(u)=K$, we have
\begin{equation*}\begin{array}{l}
\displaystyle\left\|G^{(\si)}(f)\right\|_{\al,\la}\leq\Lambda_\al(g)\sup_{t\in[0,T]}
e^{-\la t}\left(C_{\al}^{(3)} K\int_0^t\frac{(t-u)^{{
\mu}-\al}}{u^\al}du\right.\\[4mm]
\displaystyle\quad \left.+
C_{\al,T}^{(4)}\int_0^t\left((t-u)^{-2\al}+u^{-\al}\right)
 \bigg(|\si(t,u, f(u))|\right.\\[4mm]
\displaystyle\quad\quad \left.+\int_0^u\frac{|\si(t,u, f(u))-\si(t,y, f(y))|}{(u-y)^{\al+1}}dydu\bigg)\right.\\[4mm]
\displaystyle\quad
+\left.\al\int_0^t\int_0^s\int_0^u\frac{\left|\Sigma(t,s,u,y,
f)\right|}{(u-y)^{\al+1}(t-s)^{\al+1}}dyduds\right)
\\[4mm]
\displaystyle\quad \leq \Lambda_\al(g)\sup_{t\in[0,T]}\left(
C_{\al}^{(3)}KB(1-\al,{ 1+\mu}-\al) e^{-\la
t}t^{{ 1+\mu}-2\al}+ C_{\al,T}^{(4)}A_1+\al A_2\right),
\end{array}\end{equation*}
where
\begin{eqnarray*}
A_1&=&  e^{-\la t}\int_0^t\left((t-u)^{-2\al}+u^{-\al}\right)\Bigg(|\si(t,u, f(u))|\\&&\qquad\qquad+\int_0^u\frac{|\si(t,u, f(u))-\si(t,y, f(y))|}{(u-y)^{\al+1}}dydu\Bigg),\\
A_2&=&  e^{-\la
t}\int_0^t\int_0^s\int_0^u\frac{\left|\Sigma(t,s,u,y,
f)\right|}{(u-y)^{\al+1}(t-s)^{\al+1}}dyduds.
\end{eqnarray*}
Indeed, we have used (see (\ref{semprebeta})) that
\begin{equation}
\int_0^t (t-u)^q u^p
du=t^{p+q+1}\int_0^1(1-y)^{q}y^{p}dy=B(p+1,q+1)t^{p+q+1}.
\label{cotabeta}
\end{equation}
Moreover, using {\bfseries\upshape (H1)} it holds that
\begin{eqnarray} A_1&\leq&e^{-\la t}\int_0^t
((t-u)^{-2\al}+u^{-\al})\bigg(K(t^{
\mu}+u^\beta+|f(u)|)+|\si(0,0,0)|\nonumber\\&&+K\int_0^u\frac{|f(u)-f(y)|}{(u-y)^{\al+1}}dy+K\int_0^u(u-y)^{\beta-\al-1}dy\bigg)du \nonumber\\
&\leq & A_{1,1}+A_{1,2},\label{cotaA1}
\end{eqnarray}
where,
\begin{eqnarray*}
A_{1,1}&=& e^{-\la t}\int_0^t ((t-u)^{-2\al}+u^{-\al})\\&&\times \left(|\si(0,0,0)|+K\left(|f(u)|+\int_0^u\frac{|f(u)-f(y)|}{(u-y)^{\al+1}}dy\right)\right)du,\\
A_{1,2}&=& Ke^{-\la t}\int_0^t ((t-u)^{-2\al}+u^{-\al})\left(t^{\mu}
+u^\beta+\frac{1}{\beta-\al}u^{\beta-\al}\right)du.
\end{eqnarray*}
In proposition 4.2 in \cite{N-R} it has been proved that
\begin{equation} \int_0^t e^{-\la(t-u)}((t-u)^{-2\al}+u^{-\al})du
\le C_\al\la^{2\al-1} \label{sempremes}
\end{equation}
where $C_\al \le \frac{1}{1-2\alpha} +4.$ Then, the  term $A_1$ can
be treated as in \cite{N-R}. We get that
\begin{eqnarray}
\sup_{t\in[0,T]} A_{1,1}&\leq& C_\al \la^{2\al-1}\sup_{u\in[0,T]} e^{-\la u}\bigg(|\si(0,0,0)|
\nonumber\\&&+K\left(|f(u)|+\int_0^u\frac{|f(u)-f(s)|}{(u-s)^{\al+1}}ds\right)\bigg)\nonumber
\\&\leq& \la^{2\al-1} C_\al(|\si(0,0,0)|+K)(1+\left\|f\right\|_{\al,\la}).\label{cotaA11}
\end{eqnarray}
\vskip 5pt \noindent The term $A_{1,2}$ can be computed easily using
(\ref{cotabeta}). Indeed, we get that
\begin{eqnarray*}
A_{1,2}&=& K e^{-\la t}\left(\frac{t^{{
1+\mu}-2\al}}{1-2\al}+\frac{t^{{
1+\mu}-\al}}{1-\al}+\frac{t^{\beta-\al+1}}{\beta-\al+1}\right.\\[3mm]&&\left.+B(1+\beta-\al,1-2\al)
\frac{t^{\beta-3\al+1}}{\beta-\alpha}\right.\\[3mm]&&\left.+t^{\beta-2\al+1}\left(\frac{1}{(\beta-\al)(\beta-2\al+1)}+B(\beta+1,1-2\al)\right)
\right)
\\&\leq& K_{\al,\beta}^{(3)} e^{-\la t}\left(t^{
1+\mu-2\al}+t^{
1+\mu-\al}+t^{\beta-\al+1}+t^{\beta-3\al+1}+t^{\beta-2\al+1}\right),
\end{eqnarray*}
where,
\begin{eqnarray*}
K_{\al,\beta}^{(3)}&=&\frac{1}{1-2\al}+\frac{1}{1-\al}+\frac{1}{\beta-\al+1}+\frac{1}{\beta-\alpha}
B(1+\beta-\al,1-2\al)\\&&+\frac{1}{(\beta-\al)(\beta-2\al+1)}+B(\beta+1,1-2\al).
\end{eqnarray*}
So, using (\ref{sempre}) it holds that
$$\sup_{t\in[0,T]} A_{1,2} \leq K_{\al,\beta}^{(3)}K_{\al,\beta}^{(4)} ( \la^{3\al-\beta-1}{ +
\lambda^{2\alpha-\mu-1} }),$$ where
\begin{eqnarray*}
K_{\al,\beta}^{(4)}&=&\left(\frac{{
1+\mu}-2\al}{e}\right)^{{1+ \mu}-2\al}
+\left(\frac{{ 1+\mu}-\al}{e}\right)^{{
1+\mu}-\al}
+\left(\frac{\beta-\al+1}{e}\right)^{\beta-\al+1}\\&&+\left(\frac{\beta-2\al+1}{e}\right)^{\beta-2\al+1}
+\left(\frac{\beta-3\al+1}{e}\right)^{\beta-3\al+1}.
\end{eqnarray*}
Putting together (\ref{cotaA1}) and  (\ref{cotaA11}) we finally get
that
\begin{eqnarray}
\sup_{t\in[0,T]}A_1\leq \la^{2\al-1}K_{\al,\beta,K}^{(5)}\left(1+\left\|f\right\|_{\al,\la}\right), \label{cotaA1def}
\end{eqnarray}
where
\[K_{\al,\beta,K}^{(5)}=C_\al(|\si(0,0,0)|+K)+K_{\al,\beta}^{(3)}K_{\al,\beta}^{(4)}.\]

On the other hand,  we will use lemma $\ref{lema2}$ to study the term $A_2$. We can write
\begin{eqnarray*}
A_2&\leq& e^{-\la t}\int_0^t\int_0^s\int_0^u \frac{K|t-s||u-y|^\beta+ K|t-s||f(u)-f(y)|}{(u-y)^{\al+1}(t-s)^{\al+1}}dyduds\nonumber\\
&\leq& K e^{-\la t}\int_0^t\int_0^s\int_0^u |t-s|^{-\al}|u-y|^{\beta-\al-1}dyduds\nonumber\\
&&+ K e^{-\la t}\int_0^t\int_0^u\int_u^t (t-s)^{-\al}\frac{|f(u)-f(y)|}{(u-y)^{\al+1}}dsdydu\nonumber\\
&\leq& \frac{K}{\beta-\al}e^{-\la t}\int_0^t\int_0^s (t-s)^{-\al}u^{\beta-\al}duds\nonumber\\
&&+\frac{K}{1-\al}e^{-\la t}\int_0^t\int_0^u (t-u)^{1-\al}\frac{|f(u)-f(y)|}{(u-y)^{\al+1}}dydu\nonumber\\
&\leq& \frac{K}{(\beta-\al)(\beta-\al+1)}e^{-\la t}\int_0^t (t-s)^{-\al}s^{\beta-\al+1}ds\nonumber\\
&&+\frac{K}{1-\al}\left\|f\right\|_{\al,\la}\int_0^t e^{-\la (t-u)} (t-u)^{1-\al}du\nonumber\\
&\leq& \frac{K B(1-\al,\beta-\al+2)}{(\beta-\al)(\beta-\al+1)}
t^{\beta-2\al+2}e^{-\la t}+ \frac{K
T^{1-\al}}{1-\al}\la^{-1}\left\|f\right\|_{\al,\la}\nonumber
\end{eqnarray*}
and
\begin{eqnarray}
\sup_{t\in[0,T]} A_2&\leq& K_{\al,\beta}^{(6)} \la^{2\al-\beta-2}+
K_{\al,\beta,N}^{(7)}\la^{-1}\left\|f\right\|_{\al,\la}\nonumber\\&\leq&
(K_{\al,\beta}^{(6)}+K_{\al,\beta}^{(7)}) \la^{-1}
(1+\left\|f\right\|_{\al,\la}), \label{cotaA2}
\end{eqnarray}
where
\[K_{\al,\beta}^{(6)}=\frac{K B(1-\al,\beta-\al+2)}{(\beta-\al)(\beta-\al+1)}
\left(\frac{\beta-2\al}{e}\right)^{\beta-2\al+2}\quad{\textrm
and}\quad K_{\al,\beta}^{(7)}= \frac{K T^{1-\alpha}}{1-\al}.\] From
(\ref{cotaA1def}) and (\ref{cotaA2}),   (\ref{cotasi4}) holds with
\[d^{(4)}=C_{\al,K}^{(3)}K B (1-\alpha,2-\alpha) (2-2\al)^{2-2\al}e^{2\al-2}+C_{\al,T}^{(4)}K_{\al,\beta,K}^{(5)}+\al(K_{\al,\beta}^{(6)}+K_{\al,\beta}^{(7)}).\]

Assume now that $\left\|f\right\|_\infty\leq N$ and $\left\|h\right\|_\infty\leq N$. Notice that from lemma \ref{lema2} with $s_1=s_2$ we obtain that
\[|\si(t,u,f(u))-\si(t,u,h(u))-\si(s,u,f(u))+\si(s,u,h(u))|\leq K|f(u)-h(u)||t-s|.\]
Furthermore,  using (\ref{cotasi2}) with $K(u)= K |f(u)-h(u)|$, we can write,
$$\begin{array}{l}
\displaystyle\left\|G^{(\si)}(f)-G^{(\si)}(h)\right\|_{\al,\la}\\ [4mm]\displaystyle \quad\leq \Lambda_\al(g)\sup_{t\in[0,T]}e^{-\la t}\bigg(C_{\al}^{(3)} K \int_0^t\frac{|f(u)-h(u)|}{u^\al}(t-u)^{1-\al}du\\
[4mm]\displaystyle\quad\quad+C_{\al,T}^{(4)}\int_0^t ((t-u)^{-2\al}+u^{-\al})\bigg(|\si(t,u,f(u))-\si(t,u,h(u))|\\
[4mm]\displaystyle\quad\quad\quad +\int_0^u\frac{|\si(t,u,f(u))-\si(t,u,h(u))-\si(t,y,f(y))+\si(t,y,h(y))|}{(u-y)^{\al+1}}dy\bigg)du\\
[4mm]\displaystyle\quad\quad  +\al\int_0^t\int_0^s\int_0^u(u-y)^{-\al-1}(t-s)^{-\al-1}\left|\si(t,u,f(u))-\si(t,u,h(u))\right.\\
[4mm]\displaystyle\quad\quad\quad \left.-\si(s,u,f(u))+\si(s,u,h(u))-\si(t,y,f(y))+\si(t,y,h(y))\right.\\
[4mm]\displaystyle\quad\quad\quad  \left.+\si(s,y,f(y))-\si(s,y,h(y))\right|dyduds\bigg)\\
[4mm]\displaystyle\quad = \Lambda_\al(g)\sup_{t\in[0,T]}\left(C_{\al}^{(3)}K B_0+ C_{\al,T}^{(4)}(B_1+B_2)+ \alpha B_3\right),
\end{array}$$
where
\begin{eqnarray*}
B_0&=& e^{-\la t}\int_0^t\frac{|f(u)-h(u)|}{u^\al}(t-u)^{1-\al}du,\\
B_1&=& e^{-\la t}\int_0^t ((t-u)^{-2\al}+u^{-\al})|\si(t,u,f(u))-\si(t,u,h(u))|du,\\
B_2&=& e^{-\la t}\int_0^t ((t-u)^{-2\al}+u^{-\al})
\\[4mm]&&\times \int_0^u\frac{|\si(t,u,f(u))-\si(t,u,h(u))-\si(t,y,f(y))+\si(t,y,h(y))|}{(u-y)^{\al+1}}dydu,\\
B_3&=& e^{-\la t}\int_0^t\int_0^s\int_0^u(u-y)^{-\al-1}(t-s)^{-\al-1}\left|\si(t,u,f(u))-\si(t,u,h(u)) \right.\\[4mm]&&\left.-\si(s,u,f(u))+\si(s,u,h(u))-\si(t,y,f(y))+\si(t,y,h(y))\right.\\[4mm]&&\left.+\si(s,y,f(y))-\si(s,y,h(y))\right|dyduds.
\end{eqnarray*}
$B_0$ can be studied easily,
\begin{eqnarray}
\sup_{t\in[0,T]} B_0
&\leq& \left\|f-h\right\|_{\al,\la}\sup_{t\in[0,T]} t^{1-\al}\int_0^t e^{-\la (t-u)}{u^{-\al}}du\nonumber\\
&\leq& \left\|f-h\right\|_{\al,\la}\sup_{t\in[0,T]} \frac{t^{1-\al}}{\la^{1-\al}}\int_0^{\la t} e^{-x}{(\la t-x)^{-\al}}dx\nonumber\\
&\leq& \left\|f-h\right\|_{\al,\la} \frac{T^{1-\al}}{\la^{1-\al}}\sup_{z>0}\int_0^{z} e^{-x}{(z-x)^{-\al}}dx\nonumber\\
&\leq& \left\|f-h\right\|_{\al,\la} K T^{1-\al}\la^{\al-1}.
\label{cotaB0}
\end{eqnarray}
Using (\ref{sempremes}) we can deal with $B_1$
\begin{eqnarray}
\sup_{t\in[0,T]} B_1&\leq& K\sup_{t\in[0,T]}e^{-\la t}\int_0^t ((t-u)^{-2\al}+u^{-\al})|f(u)-h(u)|du\nonumber\\&\leq&
K\left\|f-h\right\|_{\al,\la}\sup_{t\in[0,T]}\int_0^t e^{-\la(t-u)}((t-u)^{-2\al}+u^{-\al})du\nonumber\\&\leq& K C_\al\la^{2\al-1}\left\|f-h\right\|_{\al,\la}, \label{cotaB1}
\end{eqnarray}
where $C_\al\leq \frac{1}{1-2\al}+4$. \vskip 5pt \noindent We will
call lemma \ref{lema1} and similar computations to those used to
deal with $B_1$ in order to obtain an estimation for $B_2$
\begin{eqnarray}
& &\sup_{t\in[0,T]} B_2 \leq \sup_{t\in[0,T]}e^{-\la t} \int_0^t
((t-u)^{-2\al}+u^{-\al})
\left(K\int_0^u\frac{\left|f(u)-h(u)\right|}{(u-y)^{\al+1-\beta}}dy\right.\nonumber\\
&&\qquad\left. + K_N\int_0^u \frac{|f(u)-h(u)|}{(u-y)^{\al+1}}\left(|f(u)-f(y)|^\delta+|h(u)-h(y)|^\delta\right)dy\right.\nonumber\\
&&\qquad \left.+ K_N\int_0^u \frac{\left|f(u)-h(u)-f(y)+h(y)\right|}{(u-y)^{\al+1}}dy\right)du\nonumber\\
&& \quad\leq \left[(K_N\left(\Delta(f)+\Delta(h)\right)+1)\int_0^t
e^{-\la(t-u)}\left((t-u)^{-2\al}+u^{-\al}\right)du\right.\nonumber
\\[4mm]&&\qquad \left.+\frac{K}{\beta-\al}\sup_{t\in[0,T]}\int_0^t e^{-\la (t-u)}((t-u)^{-2\al}+u^{-\al})u^{\beta-\al}du \right]\left\|f-h\right\|_{\al,\la}\nonumber\\
&& \quad \leq
K^{(8)}_{\al,\beta,N}\la^{2\al-1}\left\|f-h\right\|_{\al,\la}\left(1+\Delta(f)+\Delta(h)\right),\label{cotaB2}
\end{eqnarray}
where ${\displaystyle
K^{(8)}_{\al,\beta,N}=C_\al\left(K_N+\frac{KT^{\beta-\al}}{\beta-\al}\right)}$.
\vskip 5pt \noindent $B_3$ is bounded using lemma \ref{lema3}
\begin{equation} B_3 \le K_NB_{3,1}+KB_{3,2}+K_NB_{3,3},\label{cotaB3}
\end{equation}
where
\begin{eqnarray*}
B_{3,1}&=&e^{-\la t}\int_0^t\int_0^s\int_0^u \frac{|f(u)-h(u)-f(y)+h(y)|}{|t-s|^{\al}(u-y)^{\al+1}}dyduds,\\
B_{3,2}&=&e^{-\la t}\int_0^t\int_0^s\int_0^u \frac{|f(u)-h(u)|}{|t-s|^{\al}(u-y)^{\al-\beta+1}}dyduds,\\
B_{3,3}&=&e^{-\la t}\int_0^t\int_0^s\int_0^u
\frac{|f(u)-h(u)|}{|t-s|^{\al}(u-y)^{\al+1}}\left(|f(u)-f(y)|^\delta\right.\\&&\left.+\left|h(u)-h(y)\right|^\delta\right)dyduds.
\end{eqnarray*}
Then,
\begin{eqnarray}
B_{3,1}&\leq& e^{-\la t}\frac{1}{1-\al}\int_0^t\int_0^u (t-u)^{1-\al}\frac{|f(u)-h(u)-f(y)+h(y)|}{(u-y)^{\al+1}}dydu\nonumber\\
&\leq& \frac{1}{1-\al}\left\|f-h\right\|_{\al,\la}\int_0^t e^{-\la(t-u)}(t-u)^{1-\al}du\nonumber\\
&\leq& \frac{T^{1-\al}}{1-\al}\la^{-1}\left\|f-h\right\|_{\al,\la},
\label{cotaB31}
\end{eqnarray}
\begin{eqnarray}
B_{3,2}&\leq& e^{-\la t} \frac{1}{1-\al}\int_0^t\int_0^u |f(u)-h(u)|(u-y)^{\beta-\al-1}(t-u)^{1-\al}dydu\nonumber\\
&\leq&  \frac{1}{(1-\al)(\beta-\al)}\left\|f-h\right\|_{\al,\la}\int_0^t e^{-\la (t-u)} u^{\beta-\al}(t-u)^{1-\al}du\nonumber\\
&\leq&
\frac{T^{1+\beta-2\al}}{(1-\al)(\beta-\al)}\la^{-1}\left\|f-h\right\|_{\al,\la},\label{cotaB32}
\end{eqnarray}
\begin{eqnarray}
B_{3,3}&\leq&\frac{e^{-\la t}}{1-\al}\int_0^t\int_0^u\frac{|f(u)-h(u)|}{(t-u)^{\al-1}(u-y)^{\al+1}}
\nonumber\\
&  & \quad\qquad \times(|f(u)-f(y)|^\delta+|h(u)-h(y)|^\delta)dydu\nonumber\\
&\leq& \frac{1}{(1-\al)}(\Delta(f)+\Delta(h)) e^{-\la t}\int_0^t  \frac{|f(u)-h(u)|}{(t-u)^{\al-1}}du\nonumber\\
&\leq& \frac{1}{(1-\al)}(\Delta(f)+\Delta(h))  \left\|f-h\right\|_{\al,\la}   \int_0^t
\frac{e^{-\la (t-u)}}{(t-u)^{\al-1}}du\nonumber\\
&\leq&
\frac{1}{(1-\al)}\la^{-1}T^{1-\al}(\Delta(f)+\Delta(h))\left\|f-h\right\|_{\al,\la}.
\label{cotaB33}
\end{eqnarray}
Now using (\ref{cotaB3}), (\ref{cotaB31}), (\ref{cotaB32}) and (\ref{cotaB33}) we have
\begin{equation}
\sup_{t\in[0,T]} B_3\leq K^{(9)}_{\al,\beta,N}
\la^{-1}(1+\Delta(f)+\Delta(h))\left\|f-h\right\|_{\al,\lambda},
\label{cotaB3def}
\end{equation}
where
\[K^{(9)}_{\al,\beta,N}=\frac{K_N}{1-\al}+\frac{KT^{1+\beta-2\al}}{(\beta-\al)(1-\al)}+\frac{K_N T^{1-\al}}{(1-\al)}.\]
So putting together (\ref{cotaB0}), (\ref{cotaB1}), (\ref{cotaB2}) and (\ref{cotaB3def}) we obtain
that (\ref{cotasi5}) holds for
\[d_N'=\left(C_{\al,T}^{(4)}+1\right)\left(C_\al^{(3)} K K_N T^{1-\al}+KC_\al+\left(K_{\al,\beta,N}^{(8)}+K^{(9)}_{\al,\beta,N}\right)\right).\]
\hfill $\Box$
\renewcommand{\theequation}{4.\arabic{equation}}
\setcounter{equation}{0}
\section{Deterministic and stochastic equations}
Set $0<\al<\frac{1}{2}$ and $g\in W_T^{1-\al,\infty}(0,T;\R^m)$. Consider the deterministic differential equation on $\R^d$
\begin{equation}
x(t)=x_0+\int_0^t b(t,s,x(s))ds+\int_0^t \si(t,s,x(s))dg_s,\qquad t\in[0,T]. \label{eqdet}
\end{equation}
Using the notations introduced previously we can write the equation (\ref{eqdet}) as:
\begin{equation}
x(t)=x_0+F_t^{(b)}(x)+G^{(\si)}_t(x),\qquad t\in[0,T].
\end{equation}
Then we can state the result of existence and uniqueness of solution
\begin{theorem}\label{deterexi}
Assume that $b$ and $\si$ satisfy hypothesis {\bf\upshape (H1)} and {\bf\upshape (H2)} with
 $\rho=1/\al,\;\delta\leq 1 {
, \min \{ \beta, \frac{\delta}{1+\delta} \} > 1-\mu}$ and
\[{
0<1-\mu}<\al<\al_0:=\min\left\{\frac{1}{2},\beta,\frac{\delta}{1+\delta}\right\}.\]
Then  equation (\ref{eqdet}) has a unique solution $x\in
W_0(0,T;\R^d)\cap\cc^{1-\al}(0,T;\R^d)$.
\end{theorem}
The proof of this theorem follows the same computations to those
given in Theorem 5.1 in \cite{N-R}. Indeed, notice that the
estimates that we have obtained for Lebesgue integrals and
Riemann-Stieltjes integrals in Propositions \ref{Prop4.4} and
\ref{Prop4.2} are the same, with different constants,  to those
proved in  Propositions 4.4 and 4.2 of \cite{N-R}. \vskip 6pt
\noindent We provide now an upper bound for the norm of the
solution. \vskip 5pt \noindent We define $\varphi(\al,\gamma)$ as,
\[\varphi(\al,\gamma)=\qquad\begin{cases}
\frac{1}{1-\al}\qquad\mathrm{if}\; 0\leq\gamma<\frac{1-2\al}{1-\al},\\
>\frac{1}{1-2\al}\qquad\mathrm{if}\; \frac{1-2\al}{1-\al}\leq\gamma<1,\\
\frac{1}{1-2\al}\qquad\mathrm{if}\; \gamma=1.
\end{cases}\]
\begin{proposition}
Assume that $b$ and $\si$ satisfy hypothesis  of Theorem
\ref{deterexi} and that $\sigma$ satisfies also {\bf\upshape (H3)}.
Then, the unique solution of the equation (\ref{eqdet}) satisfies
\[\left\|x\right\|_{\al,\infty}\leq C_\al^{(5)} e^{C_\al^{(6)}\Lambda_\al(g)^{\varphi(\al,\gamma)}},\]
where the constants $C_\al^{(5)}$ and $C_\al^{(6)}$ depend only on $T,\,\al,\,\gamma$ and the constants that appear in conditions {\bf\upshape (H1)}, {\bf\upshape (H2)} and {\bf\upshape (H3)}.
\end{proposition}
{\bf Proof:}
We will denote by $C$ a positive constant, depending on $T,\,\al,\,\gamma$ and the constants that appear in conditions {\bf\upshape (H1)}, {\bf\upshape (H2)} and {\bf\upshape (H3)}, that will change from line to line. Using (\ref{eq1}) we have
\begin{eqnarray*}
\left|G_t^{(\si)}(x)\right|
&\leq& \Lambda_\al(g)\left(\int_0^t\frac{|\si(t,s,x(s))|}{s^\al}ds\right.\\&&+\left.\al\int_0^t\int_0^s\frac{|\si(t,s,x(s))-\si(t,r,x(r))|}{(s-r)^{\al+1}}drds\right)\\
&\leq& \Lambda_\al(g)\left(K_0\int_0^t\frac{1+|x_s|^\gamma}{s^\al}ds+\al K\int_0^t\int_0^s\frac{|x(s)-x(r)|}{(s-r)^{\al+1}}drds\right.\\&&\left.+\frac{\al K}{(\beta-\al)(\beta-\al+1)}t^{\beta-\al+1}\right)\\
&\leq& C\Lambda_\al(g)\left(1+\int_0^t\left( s^{-\al}\left|x(s)\right|^\gamma+s^{-\al}\int_0^s\frac{|x(s)-x(r)|}{(s-r)^{\al+1}}dr\right)ds\right).
\end{eqnarray*}
Then using (\ref{eq2}) for $K(u)=K$, (\ref{eq3}) and (\ref{eq4}), we
have
\begin{equation*}\begin{array}{l} \displaystyle
\int_0^t\frac{|G_t(f)-G_s(f)|}{(t-s)^{\al+1}}ds\leq\Lambda_\al(g)\left(\frac{K
B(1-\al,{ 1+\mu}-\al)}{{
\mu}-\al}t^{{
1+\mu}-2\al}\right.\\
[4mm]\displaystyle \qquad\qquad\left.+ \al\int_0^t\int_0^s\int_0^u\frac{|f(t,u)-f(s,u)-f(t,y)+f(s,y)|}{(u-y)^{\al+1}(t-s)^{\al+1}}dyduds\right.\\
[4mm]\displaystyle \qquad\qquad\left.+B(2\al,1-\al)\int_0^t\frac{|f(t,u)|}{(t-u)^{2\al}}du
\right.\\
[4mm]\displaystyle \qquad\qquad\left.+\int_0^t\int_0^u\frac{|f(t,u)-f(t,y)|}{(u-y)^{\al+1}}(t-y)^{-\al}dydu\right).
\end{array}\end{equation*}
Following the same computations we have done for the study of $A_2$ and applying the previous result combined with lemma \ref{lema2}, it holds that
\begin{equation*}\begin{array}{l}
\displaystyle
\int_0^t\frac{|G_t^{(\si)}(x)-G_s^{(\si)}(x)|}{(t-s)^{\al+1}}ds\leq
C\Lambda_\al(g)\left(
\frac{K  B(1-\al,\mu-\al+1)}{(\mu-\al)}t^{\mu-2\al+1}\right.\\
[4mm]\displaystyle \qquad\qquad \left.+
\frac{K \al B(1-\al,\beta-\al+1)}{(\beta-\al)(\beta-\al+2)}t^{\beta-2\al+2}\right.\\
[4mm]\displaystyle \qquad\qquad \left.+ \frac{K}{1-\al}\int_0^t\int_0^u (t-u)^{1-\al}\frac{|x(u)-x(y)|}{(u-y)^{\al+1}}dydu\right.\\
[4mm]\displaystyle \qquad\qquad \left.+ K_0 B(2\al,1-\al)\int_0^t\frac{1+|x(u)|^\gamma}{(t-u)^{2\al}}du\right.\\
[4mm]\displaystyle \qquad\qquad \left.+K\int_0^t\int_0^u\frac{|x(u)-x(y)|}{(u-y)^{\al+1}}(t-y)^{-\al}dydu\right.\\
[4mm]\displaystyle \qquad\qquad \left.+K\int_0^t\int_0^u(u-y)^{\beta-\al-1}(t-y)^{-\al}dydu\right)\\
[4mm]\displaystyle \qquad\leq C\Lambda_\al(g)\left(1+\int_0^t\int_0^u(t-u)^{1-\al}\frac{|x(u)-x(y)|}{(u-y)^{\al+1}}dydu\right.\\
[4mm]\displaystyle \qquad\qquad \left.+\int_0^t\frac{|x(u)|^\gamma}{(t-u)^{2\al}}du+\int_0^t\int_0^u\frac{|x(u)-x(y)|}{(u-y)^{\al+1}}(t-y)^{-\al}dydu\right)\\
[4mm]\displaystyle \qquad\leq
C\Lambda_\al(g)\left(1+\int_0^t\left(\frac{|x(u)|^\gamma}{(t-u)^{2\al}}+(t-u)^{-\al}\int_0^u\frac{|x(u)-x(y)|}{(u-y)^{\al+1}}dy\right)du\right).
\end{array}\end{equation*}
Finally using (\ref{cotamom}) we get
\begin{equation*}
\left|F^{(b)}_t(x)\right|+\int_0^t\frac{\left|F^{(b)}_t(x)-F^{(b)}_s(x)\right|}{(t-s)^{\al+1}}ds\leq C\left(1+\int_0^t(t-s)^{-\al}|x(s)|ds\right).
\end{equation*}
The proof finishes following the same computations of Proposition
5.1 in \cite{N-R}. Indeed, set
$$h(t) = \vert x(t) \vert + \int_0^t \frac{ \vert x(t) - x(s)
\vert}{ (t-s)^{\alpha + 1}} ds.$$ Then the following inequality
holds $$ h(t) \le C (1 + \Lambda_\alpha(g) ) \left( 1 + \int_0^t
\left( (t-s)^{- (1 -1/\varphi(\alpha,\gamma))}+ s^{-\alpha} \right)
h(s) ds \right).$$ Using a Gronwall inequality (see \cite{N-R}) we
finish the proof. \hfill $\Box$ \vskip 7pt

\noindent We will finish applying the results for deterministic equations to stochastic equations.
The stochastic integral appearing throughout this paper is a path-wise Riemann-Stieltjes integral and it is
well know that this integral exists if the process that we integrate has
H\"older continuous trajectories of order larger than $1-H$.
\\
Set $\al \in (1-H, \frac12)$. For any $\delta \in (0,2)$, by
Fernique's theorem it holds that $$ E( \exp(\Lambda_\al (W)^\delta
)) < \infty.$$
Then if $u_t=\{u_t(s), s \in [0,T]\}$ is a stochastic
process whose trajectories belong to the space $W_T^{\al,1}
(0,T)$, the Riemann-Sieltjes integral $\int_0^T u_t(s) dW_s$ exists
and we have that
$$ \left| \int_0^T u_t(s) dW_s\right| \le G \Vert u_t \Vert_{\al,1},$$
where $G$ is a random variable with moments of all orders (see Lemma
7.5 in \cite{N-R}). Moreover, if the trajectories of $u_t$ belong to
$W_0^{\al,\infty} (0,T)$, then the indefinite integral $\int_0^T
u_t(s) dW_s $ is H\"older continuous of order $1-\al$ and with
trajectories in $W_0^{\al,\infty} (0,T)$. As a simple consequence of
these facts and our deterministic results, we obtain the proof of
our main theorem.

\renewcommand{\theequation}{5.\arabic{equation}}
\setcounter{equation}{0}
\section{Appendix}

In this section we give some technical lemmas used in the estimates
obtained in Section 3.

\begin{lemma} \label{lema1}
Let $\si:[0,T]^2\times\R\rightarrow\R$ be a function satisfying the hypothesis {\bf\upshape (H1)}. Then for all $N>0$, $t,s_1,s_2\in\R$ and $|x_1|,|x_2|,|x_3|,|x_4|\leq N$
\begin{equation}\begin{array}{l}
\displaystyle
|\si(t,s_1,x_1)-\si(t,s_2,x_2)-\si(t,s_1,x_3)+\si(t,s_2,x_4)|\\
[4mm]\displaystyle \qquad\qquad\qquad\leq
K_N|x_1-x_2-x_3+x_4|+K|x_1-x_3||s_2-s_1|^\beta\\ [4mm]\displaystyle
\qquad\qquad\qquad\quad+K_N|x_1-x_3|\left(|x_1-x_2|^\delta+
|x_3-x_4|^\delta\right).
\end{array}\label{le1}\end{equation}
\end{lemma}
{\bf Proof:}
By the mean value theorem we can write
\begin{equation*}\begin{array}{l}
\displaystyle
\si(t,s_1,x_1)-\si(t,s_2,x_2)-\si(t,s_1,x_3)+\si(t,s_2,x_4)\\ [4mm]\displaystyle \qquad= \int_0^1 (x_1-x_3)\partial_x \si(t,s_1,\theta x_1+ (1-\theta) x_3)d\theta\\ [4mm]\displaystyle \qquad\quad-\int_0^1 (x_2-x_4)\partial_x \si(t,s_2,\theta x_2+ (1-\theta) x_4)d\theta\\ [4mm]\displaystyle \qquad =\int_0^1 (x_1-x_2-x_3+x_4)\partial_x\si(t,s_2,\theta x_2+(1-\theta)x_4)d\theta\\ [4mm]\displaystyle \qquad\quad + \int_0^1 (x_1-x_3)\left(\partial_x(t,s_1,\theta x_1+(1-\theta)x_3)\right.\\ [4mm]\displaystyle \qquad\qquad \left.-\partial_x\si(t,s_2,\theta x_2+(1-\theta)x_4)\right)d\theta.
\end{array}\end{equation*}
We can obtain (\ref{le1}) using the hypothesis {\bf{(H1)}}. \hfill
$\Box$
\begin{lemma} \label{lema2}
Let $\si:[0,T]^2\times\R\rightarrow\R$ be a function satisfying the hypothesis {\bf\upshape (H1)}. Then for all $t_1,t_2,s_1,s_2\in\R$ and 
$x_1,x_2\in \R$
\begin{equation}\begin{array}{l}
\displaystyle
|\si(t_1,s_1,x_1)-\si(t_2,s_1,x_1)-\si(t_1,s_2,x_2)+\si(t_2,s_2,x_2)|\\ [4mm]\displaystyle \qquad\qquad\quad\leq\; K|t_1-t_2|\left(|s_1-s_2|^\beta+|x_1-x_2|\right).
\end{array}\label{le2}\end{equation}
\end{lemma}
{\bf Proof:}
By the mean value theorem we can write
\begin{equation*}\begin{array}{l}
\displaystyle
\si(t_1,s_1,x_1)-\si(t_2,s_1,x_1)-\si(t_1,s_2,x_2)+\si(t_2,s_2,x_2)\\ [4mm]\displaystyle \qquad= \int_0^1
(t_1-t_2)\partial_t \si(\theta t_1+ (1-\theta) t_2,s_1,x_1)d\theta\\ [4mm]\displaystyle
 \qquad\quad-\int_0^1 (t_1-t_2)\partial_t \si(\theta t_1+ (1-\theta) t_2,s_2,x_2)d\theta\\ [4mm]\displaystyle
  \qquad =\int_0^1 (t_1-t_2)\left(\partial_t\si(\theta t_1+ (1-\theta) t_2,s_1,x_1)\right.\\ [4mm]\displaystyle
  \qquad\qquad\left.-\partial_t\si(\theta t_1+ (1-\theta) t_2,s_2,x_1)\right)d\theta\\ [4mm]\displaystyle
   \qquad\quad + \int_0^1 (t_1-t_2)\left(\partial_t(\theta t_1+ (1-\theta) t_2,s_2,x_1)\right.\\
   [4mm]\displaystyle \qquad\qquad\left.-\partial_t\si(\theta t_1+ (1-\theta) t_2,s_2,
   x_2)\right)d\theta,
\end{array}\end{equation*}
and we can obtain (\ref{le2}) using the hypothesis {\bf{(H1)}}.
\hfill $\Box$
\begin{lemma} \label{lema3}
Let $\si:[0,T]^2\times\R\rightarrow\R$ be a function satisfying the hypothesis {\bf\upshape (H1)}. Then for all $N>0$, $t_1,t_2,s_1,s_2\in\R$ and $|x_1|,|x_2|,|x_3|,|x_4|\leq N$
\begin{equation}\begin{array}{l}
\displaystyle
|\si(t_1,s_1,x_1)-\si(t_1,s_1,x_2)-\si(t_2,s_1,x_1)+\si(t_2,s_1,x_2)\\
[4mm]\displaystyle \qquad\qquad\quad\,
-\si(t_1,s_2,x_3)+\si(t_1,s_2,x_4)+\si(t_2,s_2,x_3)-\si(t_2,s_2,x_4)|\\
[4mm]\displaystyle \qquad\quad\quad\leq
K_N|t_1-t_2||x_1-x_2-x_3+x_4|+K|x_1-x_2||t_1-t_2||s_1-s_2|^\beta \\
[4mm]\displaystyle \qquad\qquad\qquad\ +
K_N|x_1-x_2||t_1-t_2|\left(|x_1-x_3|^\delta+|x_2-x_4|^\delta\right).
\end{array}\label{le3}\end{equation}
\end{lemma}
{\bf Proof:}
By the mean value theorem we can write
\begin{equation*}\begin{array}{l}
\displaystyle
\si(t_1,s_1,x_1)-\si(t_1,s_1,x_2)-\si(t_2,s_1,x_1)+\si(t_2,s_1,x_2)-\si(t_1,s_2,x_3)\\
[4mm]\displaystyle \qquad\quad\qquad
+\si(t_1,s_2,x_4)+\si(t_2,s_2,x_3)-\si(t_2,s_2,x_4)\\
[4mm]\displaystyle \qquad
= (x_1-x_2)\int_0^1 \left[\partial_x \si(t_1,s_1,\theta x_1+(1-\theta)x_2)\right.\\[4mm]\displaystyle \qquad\qquad\left.-\partial_x \si(t_2,s_1,\theta x_1+ (1-\theta)x_2)\right]d\theta\\ [4mm]\displaystyle \qquad\quad
-(x_3-x_4)\int_0^1 \left[\partial_x \si(t_1,s_2,\theta x_3+(1-\theta)x_4)\right.\\[4mm]\displaystyle \qquad\qquad\left.-\partial_x \si(t_2,s_2,\theta x_3+ (1-\theta)x_4)\right]d\theta\\ [4mm]\displaystyle \qquad
=(t_1-t_2)\bigg[(x_1-x_2)\\ [4mm]\displaystyle
\qquad\qquad\int_0^1\int_0^1  \partial_{x,t}^2\si(\la t_1+ (1-\la)
t_2,s_1,\theta x_1+(1-\theta)x_2)d\la d\theta\\ [4mm]\displaystyle
\qquad\quad - (x_3-x_4)\int_0^1\int_0^1  \partial_{x,t}^2(\la t_1+
(1-\la) t_2,s_2,\theta x_3+(1-\theta)x_4)d\la d\theta\bigg]\\
[4mm]\displaystyle \qquad \leq
(t_1-t_2)\left[(x_1-x_2)\int_0^1\int_0^1\left(\partial_{x,t}^2(\la
t_1+ (1-\la) t_2,s_1,\theta x_1+(1-\theta)x_2)\right.\right.\\
[4mm]\displaystyle \qquad\qquad \left.\left.-\partial_{x,t}^2(\la
t_1+ (1-\la) t_2,s_2,\theta x_3+(1-\theta)x_4)\right)d\la
d\theta\right.\\ [4mm]\displaystyle \qquad\quad \left.+
(x_1-x_2-x_3+x_4)\right.\\ [4mm]\displaystyle
\qquad\qquad\left.\int_0^1\int_0^1\partial_{x,t}^2\si(\la t_1+
(1-\la) t_2,s_2,\theta x_3+(1-\theta)x_4)d\la d\theta\right].
\end{array}\end{equation*}
So we can obtain (\ref{le3}) using the hypothesis {\bf{(H1)}}.
\hfill $\Box$

\section*{Acknowledgements}
This work was partially supported   by DGES Grant MTM09-07203
(Mireia Besal\'u and Carles Rovira).

\end{document}